\documentclass[12pt]{article}

\usepackage{sectsty}
\usepackage{fullpage}
\usepackage[export]{adjustbox}
\usepackage{graphicx}

\usepackage{cancel,enumerate}
\usepackage[affil-it]{authblk}

\usepackage{amsmath}
\usepackage{amstext, amssymb}
\usepackage[mathscr]{euscript}

\newcommand{\diint}{\displaystyle\int}

\DeclareMathOperator\arctanh{arctanh}
\usepackage{mathtools}

\title{Fisher Transformation via Edgeworth Expansion}
\author{Jan Vrbik}
\affil{Department of Mathematics and Statistics\\ Brock University, Canada}
\date{\today}

\begin{document}

\maketitle
\begin{abstract}
We show how to calculate individual terms of the Edgeworth series to
approximate the distribution of the Pearson correlation coefficient with the
help of a simple Mathematica program. We also demonstrate how to eliminate
the corresponding skewness, thus making the approximation substantially more
accurate. This leads, in a rather natural way, to deriving a superior (in
terms of its accuracy) version of Fisher's $z$ transformation. The code can
be easily modified to deal with any sample statistics defined as a function
of several sample means, based on a random independent sample from a
multivariate distribution.

Keywords: Edgeworth series, Pearson correlation, Fisher transformation,
moment generating function
\end{abstract}

\section{Introduction}

This article is an extension of the original work of \cite{fisher}, 
\cite{hotel} and \cite{winter} (the last one correctly discounting the importance of `variance stabilization'), with the aim of explicitly deriving \emph{all} terms
necessary to achieve an $O(n^{-3/2})$ accuracy of the desired approximation. The
resulting formula then becomes substantially more accurate than those found in the
existing literature.

One of the most natural ways to derive the Fisher transformation of the
empirical correlation coefficient $r$ (assuming a sample from a bivariate
Normal distribution) is to start by finding the first few terms of the Edgeworth
expansion \cite{edge} of the sampling distribution of an arbitrary function
of $r$ say $G(r)$. This is done by standardizing $G(r)$ by an $O(n^{-2})$%
-accurate (in terms of its error) expected value $m$ and $O(n^{-3})$%
-accurate variance $V$, and finding its $O(n^{-3/2})$-accurate skewness $%
\Gamma_{3}$ and $O(n^{-2})$-accurate \emph{excess} (subtracting $3$)
kurtosis $\Gamma_{4}$ (where $n$ is the sample size). The approximate
probability density function (PDF) of%
\begin{equation}
Z\coloneqq\frac{G(r)-m}{\sqrt{V}}  \label{stand}
\end{equation}
is then given by
\begin{equation}
f_{Z}(z) \coloneqq 
\frac{\exp \left( -\dfrac{z^{2}}{2}\right)}{\sqrt{2\pi }}
\times \left(1+\frac{\Gamma_{3}(z^{3}-3z)}{6}+\frac{\Gamma_{4}(z^{4}-6z^{2}+3)}{24}+\frac{\Gamma_{3}^{2}(z^{6}-15z^{4}+45z^{2}-15)}{72}\right)  \label{edge}
\end{equation}%
and has an $O(n^{-3/2})$-proportionate error (compared to an $O(n^{-1/2})$
error of the basic Normal approximation).

It is then possible to find $G$ (based on the resulting differential
equation) to make $\Gamma_{3}$ equal to zero (to the $O(n^{-3/2}$) level of
accuracy) for any value of the `true' correlation coefficient $\rho$. This
yields the expected `$\arctanh$' transformation, but also suggests a subtle
correction to it, making the resulting approximation substantially more
accurate.

\section{Key concepts and formulas}

In this section we assume sampling from \emph{any} specific multivariate
distribution; to simplify the notation, our definitions and examples will
use a tri-variate case and call the random variables $X_{1}$, $X_{2}$ and $%
X_{3}$, (generalizing is easy). The corresponding \emph{central} \textsc{%
moments} are then defined by%
\begin{equation}
\mu_{i,j,k} \coloneqq \mathbb{E}\left[ (X_{1}-\mu_{1})^{i}(X_{2}-\mu_{2})^{j}(X_{3}-\mu_{3})^{k}\right]  \label{moment}
\end{equation}%
where $i+j+k$ is the moment's \textsc{order}. They can be conveniently
computed based on the moment generating function (MGF) of the corresponding 
\emph{centralized} random variables, namely from%
\begin{align}
M(t_{1},t_{2},t_{3}) 
&\coloneqq \mathbb{E}\left[ \exp (t_{1}(X_{1}-\mu_{1})+t_{2}(X_{2}-\mu_{2})+t_{3}(X_{3}-\mu_{3}))\right]\label{MGF} \\
&\equiv 1+\sum_{i+j+k\geq 2}^{\infty }\frac{\mu_{i,j,k}}{i!\, j!\, k!} t_{1}^{i}t_{2}^{j}t_{3}^{k}  \notag
\end{align}
by differentiating it with respect to $t_{1}$, $t_{2}$ and $t_{3}$
correspondingly $i$, $j$ and $k$ times, and then setting each $t_{\ell }$
equal to zero. If an explicit formula for such an MGF cannot be found (the
integration may have no analytic answer), it is sufficient to replace it by
the last line of (\ref{MGF}), with the summation truncated to exclude terms
beyond the fourth order (computationally more feasible).

It is well known (and easy to prove) that, when replacing the individual
random variables by their respective \emph{sample} means, the MGF of the new
set, namely of $\overline{X_{1}}-\mu_{1}$, $\overline{X_{2}}-\mu_{2}$, 
$\overline{X_{3}}-\mu_{3}$, is given by
\begin{equation}
M\left( \frac{t_{1}}{n},\frac{t_{2}}{n},\frac{t_{3}}{n}\right)^{n}
\label{MGFM}
\end{equation}%
This implies that finding%
\begin{equation}
\bar{\mu}_{i,j,k}\coloneqq\mathbb{E}\left[ (\overline{X_{1}}-\mu_{1})^{i}(\overline{X_{2}}-\mu_{2})^{j}(\overline{X_{3}}-\mu_{3})^{k}\right]
\label{momM}
\end{equation}%
can be achieved by the same $i$, $j$, $k$-fold differentiation of (\ref{MGFM}%
) and subsequent $t_{\ell }=0$ substitution as before. Thus, for example, we
get%
\begin{equation}
\bar{\mu}_{2,1,1}=\frac{\mu_{2,0,0}\mu_{0,1,1}+2\mu_{1,1,0}\mu_{1,0,1}}{n^{2}}+\frac{\mu_{2,1,1}-\mu_{2,0,0}\mu_{0,1,1}-2\mu_{1,1,0}\mu_{1,0,1}}{n^{3}}
\end{equation}%
etc. Not surprisingly, the complexity of these formulas increases
`exponentially' with the moment's order.

To find an $O(n^{-3/2})$-accurate approximation to the PDF of any function
of sample means, say 
$H\left(\overline{X_{1}},\overline{X_{2}},\overline{X_{3}}\right)$,
we need to re-write this function as%
\begin{equation}
H\left(\mu_{1}+\varepsilon (\overline{X_{1}}-\mu_{1}),\mu_{2}+\varepsilon (\overline{X_{2}}-\mu_{2}),\mu_{3}+\varepsilon (\overline{X_{3}}-\mu_{3})\right)  \label{eps}
\end{equation}%
and expand it in $\varepsilon $, up to and including $\varepsilon^{3}$
terms (with the understanding that $\varepsilon $ will be set to $1$
eventually); note that the first (constant) term of this expansion is 
$H\left( \mu_{1},\mu_{2},\mu_{3}\right) $.

We then compute the expected value of the first \emph{four} powers of the
result, after subtracting $H\left( \mu_{1},\mu_{2},\mu_{3}\right) $ from
it (to simplify the corresponding algebra). This requires further expanding
of these powers in $\varepsilon $, up to and including $\varepsilon^{2}$, 
$\varepsilon^{4}$, $\varepsilon^{4}$ and $\varepsilon^{6}$ terms
respectively, before applying (\ref{momM}) to the individual terms.

The resulting four moments of 
$H\left( \overline{X_{1}},\overline{X_{2}},\overline{X_{3}}\right) -H\left( \mu_{1},\mu_{2},\mu_{3}\right)$
are then easily converted to the corresponding $m$, $V$, $\Gamma_{3}$ and 
$\Gamma_{4}$; these should be further simplified by keeping only the leading
terms of their $\frac{1}{n}$ expansion (with the exception of $V$, where
both $\frac{1}{n}$ and $\frac{1}{n^{2}}$-proportional terms are needed).
This ensures that all terms contributing to the final $O(n^{-3/2})$ accuracy
of the final answer are included, while the rest of them (often incorrect, since
incomplete) have been eliminated.

\section{Examples}

Pearson's correlation coefficient $r$ is defined, using the $\varepsilon$
notation of the previous section, as
\begin{equation}
r
\coloneqq
\frac
{\rho +\varepsilon (\overline{XY}-\rho )-\varepsilon^{2}\overline{X}\cdot \overline{Y}}
{\sqrt{\left( 1+\varepsilon (\overline{X^{2}}-1)-\varepsilon^{2}\overline{X}^{2}\right) \left( 1+\varepsilon (\overline{Y^{2}}-1)-\varepsilon^{2}\overline{Y}^{2}\right) }}  \label{r}
\end{equation}%
Note that $r$ is a function of \emph{five} different sample means. To
investigate its sampling distribution, we assume a random independent sample
of size $n$ from a bivariate Normal distribution of $X$ and $Y$, with both
means equal to $0$ and both variances equal to $1$, while the true (or
`population') correlation coefficient is $\rho $ (the $r$ distribution is
the same whatever means and variances we use; we have thus made the simplest
choice).

\begin{figure}[htbp]
\begin{center}
\includegraphics[width=0.8\textwidth]{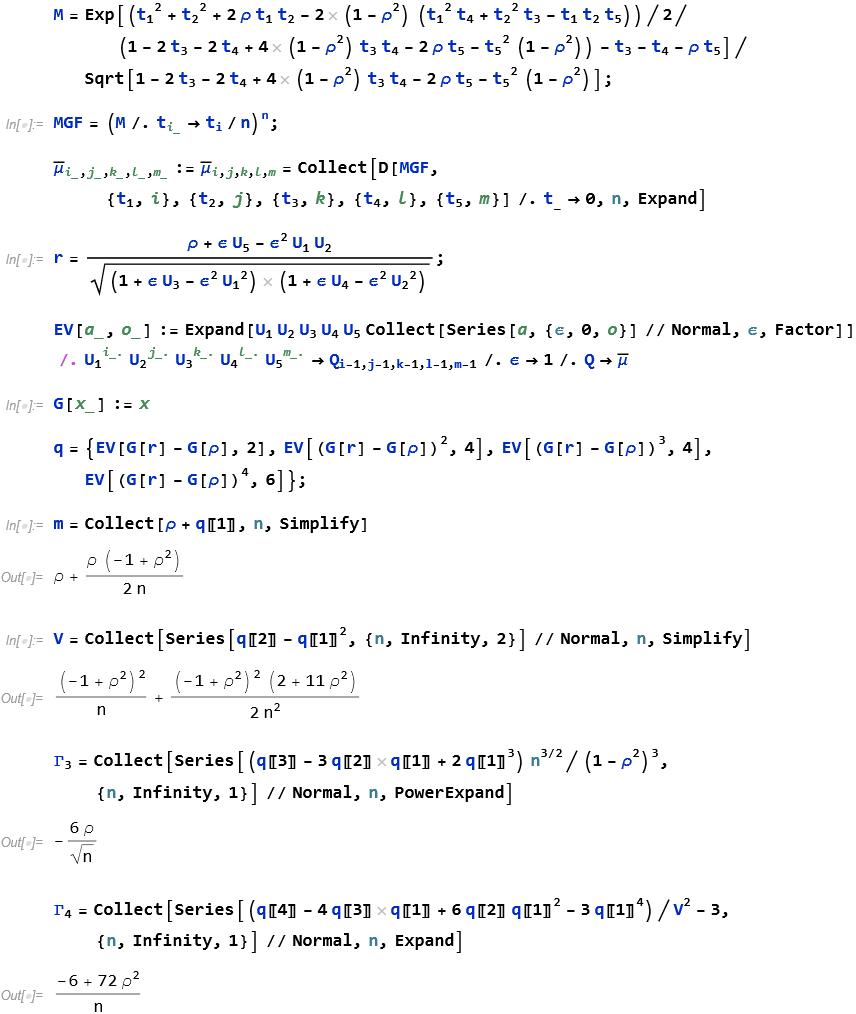}
\caption{Mathematica code}
\label{fig::code}
\end{center}
\end{figure}

The MGF of the required five centralized variables, namely $X$, $Y$, $X^{2}-1$, 
$Y^{2}-1$ and $XY-\rho$, is the result of the following (rather
routine) double integration
\begin{equation}
\frac{\diint\limits_{-\infty }^{\infty }\exp \left( -\dfrac{x^{2}+y^{2}-2\rho xy}{2(1-\rho^{2})}+t_{1}x+t_{2}y+t_{3}(x^{2}-1)+t_{4}(y^{2}-1)+t_{5}(xy-\rho )\right) }{2\pi \sqrt{1-\rho^{2}}}
\end{equation}%
quoted in the Mathematica code of Figure 1 (which calls it $M$). This
expression is then easily converted into the MGF of $\overline{X}$, 
$\overline{Y}$, $\overline{X^{2}}-1$, $\overline{Y^{2}}-1$ and 
$\overline{XY}-\rho$, and the first four moments of $r-\rho $ (no further transformation
is applied to $r$ in this example) are found; the corresponding $m$, $V$, 
$\Gamma_{3}$ and $\Gamma_{4}$, properly truncated in their $\frac{1}{n}$
expansions, then easily follow. The complete program and its output are
displayed in Figure \ref{fig::code}.

\begin{figure}[htbp]
\begin{center}
\includegraphics[width=0.8\textwidth]{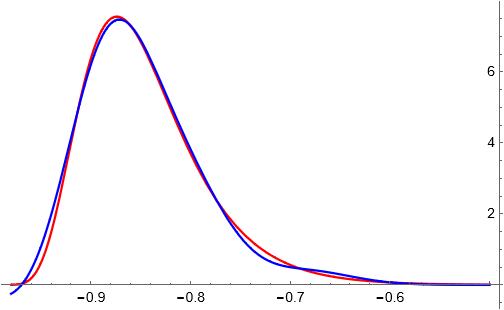}
\caption{Exact (red) and Edgeworth-series (blue) PDF of $r$}
\label{fig1}
\end{center}
\end{figure}

The resulting approximation is than given by (\ref{stand}) and (\ref{edge});
the PDF of $Z$ can be easily converted to an approximate PDF of $r$ by
\begin{equation}
f\left( r\right) \simeq \frac{f_{Z}\left( \dfrac{r-m}{\sqrt{V}}\right) }{%
\sqrt{V}}
\end{equation}%
Using $n=35$ and $\rho =-0.85$, we show (in Figure \ref{fig1} - the exact PDF is red,
the approximate one is blue) how the approximation compares to the exact
answer (which, rather atypically for a sample statistic of this complexity,
has an analytic form - see \cite{hotel}). The maximum error of this approximation, realized when computing $\Pr(-0.9685<r<-0.9133)$, is less than $0.8\%$.

\subsection{Fisher transformation}

Re-running the same program with the `$G[x\_]\coloneqq x$' line removed yields, for $%
\Gamma_{3}$, the following expression%
\begin{equation}
\frac{3G^{\prime }[\rho ]\left( (1-\rho^{2})G^{\prime \prime }[\rho ]-2\rho
~G^{\prime }[\rho ]\right) }{\sqrt{n}}
\end{equation}%
To make it zero (for any value of $\rho $), $G$ needs to be either a
constant (which is clearly inadmissible) or a solution to $(1-\rho
^{2})G^{\prime \prime }(\rho )-2\rho ~G^{\prime }(\rho )=0$, namely 
$G(\rho)=\arctanh(\rho )$, having chosen its simplest form (all other
possibilities would yield the same $Z$).

Running the same program one more time with `$G[x\_]\coloneqq$~ArcTanh$[x]$' then
results in 
\begin{align*}
m &= \arctanh(\rho )+\frac{\rho }{2n} \\
V &= \frac{1}{n}+\frac{6-\rho^{2}}{2n^{2}} \\
\Gamma_{3} &= 0\ \ \ \ \ \ \text{(by design)} \\
\Gamma_{4} &= \frac{2}{n}
\end{align*}%
This time%
\begin{equation}
f\left( r\right) \simeq \frac{f_{Z}\left( \dfrac{\arctanh(r)-m}{\sqrt{V}}\right) }{(1-r^{2})\sqrt{V}}
\end{equation}%
which, when using the previous choice of $n$ and $\rho $, yields a PDF
visually indistinguishable from the exact answer; its maximum error (now
involving a large interval of values - thus not the fairest of comparisons)
is $0.36\%$. This accuracy is maintained (actually, rather fortuitously
reduced to $0.05\%$) after dropping the $\Gamma_{4}$ term, whose
contribution is (for this transformation of $r$) practically negligible, due
to its relatively small size.

\begin{figure}[htbp]
\begin{center}
\includegraphics[width=0.8\textwidth]{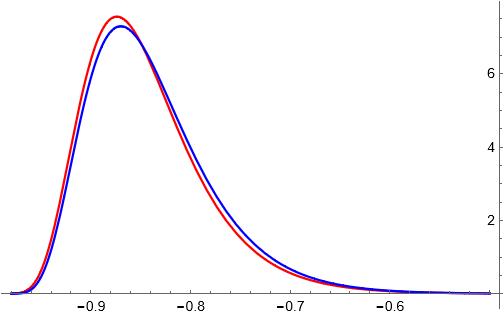}
\caption{Exact (red) and basic-Fisher (blue) PDF of $r$}
\label{fig::fig2}
\end{center}
\end{figure}

Note that the `basic' Fisher transformation similarly ignores $\Gamma_{4}$.
but it also drops the $\frac{1}{n}$-proportional correction to $m$ and uses $%
V=\dfrac{1}{n-3}$ in place of our result. This affects, quite adversely, its
accuracy, as seen in Figure \ref{fig::fig2}.%

The largest error of this approximation is over $3.5\%$ - clearly
inacceptable!

These two examples should suffice to illustrate using the Edgeworth series
in situations going well beyond its original formulation and intended
purpose.

\section{Conclusion}

We have delineated a procedure for constructing an accurate approximation
for a PDF of any function of several sample means, when sampling a specific
univariate or multivariate distribution. It is based on finding, to a
specific accuracy (in terms of their $\frac{1}{n}$ expansions), the
corresponding mean, variance, skewness and excess kurtosis; in the case of a
single-parameter distribution, it is usually possible to find a
transformation of the sample statistics which eliminates skewness, thus
making the approximation both simpler and more accurate. We have tested the
technique against the `classical' example of Fisher transformation
(suggesting a minor modification leading to a significant improvement), but
its main applicability is to situations with no exact solution (constructing
an \emph{approximate} PDF is then the best we can do). We should mention that we
have not attempted to optimize the algorithm (this would require introducing
cumulants), but since our program takes only a few seconds to execute, this
would not appear necessary in most cases of interest.

\end{document}